\documentclass[12pt]{article}
\usepackage{amsmath}
\usepackage{amssymb}
\usepackage{url}

\usepackage[T1]{fontenc}
\usepackage{mathptmx}
\usepackage{diagrams}
\diagramstyle{midshaft,PostScript=dvips,nohug}
\newarrow{Dotsto}.. ..>

\newtheorem{lemma}{Lemma}[section]
\newtheorem{prop}[lemma]{Proposition}

\newtheorem{defn}[lemma]{Definition}

\newcommand{\K}{\mathcal{K}}
\newcommand{\X}{\mathcal{X}}

\newcommand{\B }{\mathcal{B}}

\newcommand{\dcoc}{\triangleleft}

\title{The dual fibration in elementary terms}
\author{Anders Kock}
\date{}

\begin{document}

\maketitle

We give an elementary construction of the dual fibration of a 
fibration. It does not use the non-elementary notion of (pseudo-) functor into the category of 
categories. In fact, it is clear that the construction we present 
makes sense for internal categories and fibrations in any exact 
category.

The dual fibration of a fibration $\X \to \B$ over $\B$ is described in e.g.\ [Borceux] 
II.8.3 via a pseudofunctor $F: \B ^{op} \to Cat$ (the category of 
categories), by composing $F$ with the (covarariant!) dualization 
functor $Cat \to Cat$; choosing such an $F$ is tantamount to choosing 
a cleavage for the fibration. In the present section, we give  an 
alternative description of the dual fibration, which is 
elementary 
and choice-free.

\section{Fibrations}

We recall here some classical notions. 

Let $\pi : \X \to \B$ be any functor. For $\alpha :A \to B$ in $\B$, 
and for objects $X, Y\in \X$ with $\pi (X)=A$ and $\pi (Y)=B$,
let $\hom_{\alpha}(X,Y)$ be the set of arrows $h:X\to Y$ in $\X$ with $\pi 
(h)=\alpha$.  
For any arrow $\xi : C\to A$, and any object $Z\in \X$ 
with $\pi (Z)=C$, post-composition with $h$ defines a map $$h_{*}:\hom 
_{\xi}(Z,X)\to \hom_{\xi . \alpha} (Z,Y).$$
(we compose from left to right). Recall that $h$ is called {\em 
Cartesian} if this map is a bijection, for all such $\xi$ and $Z$.

If $h$ is Cartesian, the injectivity of $h_{*}$ implies the 
cancellation property that $h$ is ``monic w.r.\ to $\pi$'', 
meaning that  for parallel arrows $k,k'$ in $\X$ with 
codomain $X$, and with $\pi (k)=\pi (k')$, we have that 
$k.h = k'.h \mbox { implies } k=k'$. 

For later use, we recall a basic fact: 
\begin{lemma}\label{lemmax} If $k=k'.h$ is Cartesian, and 
$h$ is Cartesian then $k'$ is Cartesian.\end{lemma}

The functor $\pi :\X \to \B$ is called a {\em fibration} if for every 
$\alpha : A \to B$ in $\B$ and any $Y \in \X$ with $\pi (Y)=B$, there 
exists a Cartesian arrow over $\alpha$ with codomain $Y$. The {\em 
fibre} over $A\in \B$ is the category whose objects are the $X \in 
\X$ with $\pi (X)=A$, and whose arrows are arrows in $\X$ 
which by $\pi$ map to  $1_{A}$; such arrows are 
called {\em vertical} (over $A$).

\medskip

All this is standard, dating back essentially to early French 
category theory (Grothendieck, Chevalley, Giraud, B\'{e}nabou,\ldots).
For a modern account, see \cite{Borceux} II.8.1,  \cite{Johnstone}
 B.1.3, or 
\cite{Kock}. 
Note that these notions are ele\-mentary (they make sense for category 
objects in any left exact category), and they do not depend on the 
non-elementary notions of {\em cleavage}, or $Cat$-valued {\em pseudofunctor}.

\section{The ``factorization system'' for a fibration}

In the diagrams  below, we try to make display vertical arrows 
vertically, and Cartesian arrows horizontally. 

   Recall from the literature that if $\pi: \X \to \B$ is a fibration, then every arrow $z$ in $\X$ may be written as a composite 
of a vertical arrow followed by a cartesian arrow. And, crucially, 
this decomposition of $z$ is unique modulo a {\em unique} vertical
isomorphism. Or, equivalently, modulo an arrow which is at the same 
time vertical and cartesian. (Recall that for vertical arrows, 
cartesian is equivalent to isomorphism (= invertible).) This means that every arrow $z$ in $\X$ 
may be represented by a pair $(v,h)$ of arrows with $v$ vertical and 
$h$ cartesian, with $z=v.h$. Thus the codomain of $v$ is the domain 
of $h$. We call such a pair a ``vh composition pair'', to make the 
analogy with vh spans, to be considered below, more explicit. Two such pairs $(v,h)$ and 
$(v',h')$ represent the same arrow iff there exists a vertical 
cartesian (necessarily unique, and necessarily invertible) $i$ such that 
\begin{equation}\label{eqx1}v.i=v' \mbox{ and } i.h'=h. \end{equation}
We say that $(v,h)$ and $(v',h')$ are {\em equivalent} if this holds.
The composition of arrows in $\X$ can be described in terms of 
representative vh composition pairs, as follows. If $z_{j}$ is represented by 
$(v_{j},h_{j})$ for $j=1, 2$, then $z_{1}.z_{2}$ is represented by 
$(v_{1}.w, k.h_{2})$, where $k$ is cartesian over $\pi (h_{1})$ and $w$ 
is vertical, and the square displayed commutes:
$$\begin{diagram}\cdot&&&&&\\
\dTo^{v_{1}}&&&&\\
\cdot & \rTo ^{h_{1}}&\cdot&&&\\
\dTo^{w}&&\dTo_{v_{2}}&&\\
\cdot&\rTo_{k}&\cdot & \rTo_{h_{2}}&\cdot
\end{diagram}$$
Such $k$ and $w$ exists (uniquely, up to unique vertical cartesian 
arrows): 
construct first $k$ as a cartesian lift of $\pi (h_{1})$, then use 
the universal property of cartesian arrows to construct $w$.

The arrows $z_{1}$ and $z_{2}$ may be inserted, completing the 
diagram with two commutative triangles, since $z_{j}=v_{j}.h_{j}$. 
But if we refrain from doing so, we have a blueprint for 
 a succinct 
and choice-free description of the fibrewise dual $\X ^{*}$ of the fibration 
$\X\to \B$. 

Note that a vh  factorization of an arrow in $\X$ is much reminiscent 
of the  factorization for an $E$-$M$ factorization system, as in 
[Borceux] I.5.5, say, (with the 
class of  
vertical arrows playing the role of $E$, and the class of cartesian 
arrows  playing the 
role of $M$; however, note that not every isomorphism in $\X$ is 
vertical.

\section{The dual fibration $\X^{*}$}

The construction presented in this Section is still elementary, but 
requires more than just left exactness in the category where it is 
performed, namely exactness; this implies that good quotients exist 
for equivalence relations, and that maps on such a quotient can be 
defined by assigning values on representative elements for the 
equivalence classes. -- We present the construction in the exact 
category of sets, for simplicity.

Given a fibration $\pi: \X \to \B$. We describe another category 
$\X^{*}$ over $\B$, as follows: The objects of $\X^{*}$ are the same as those of $\X$; 
the arrows $X \to Y$ are represented by vh spans, in the following 
sense: 
\begin{defn}A {\em vh span} in $\X$ from $X$ to $Y$ is a diagram in $\X$ of the form
\begin{equation}\label{vhx}\begin{diagram}\cdot & \rTo^{h}& Y\\
\dTo^{v}&&\\
X&&
\end{diagram}\end{equation}
with $v$ vertical and $h$  cartesian.
\end{defn}
The set of arrows in $\X^{*}$ from $X$ 
to $Y$ are equivalence classes of vh spans from $X$ to $Y$, for the 
equivalence relation $\equiv$ given by
$(v,h) \equiv (v',h')$ if there exists a vertical isomorphism $i$ 
(necessarily unique) in 
$\X$ so that 
\begin{equation}\label{eqx2}i.v.=v' \mbox{ and } i.h=h'. \end{equation}
We denote the equivalence class of the vh span $(v,h)$ by $\{(v,h)\}$. 
They are the arrows of $\X^{*}$;   the direction of a the arrow $\{(v,h)\}$ 
is determined by its 
cartesian part $h$.
  
 Composition  has to be described in terms of  representative 
pairs; it is in fact the standard composite of spans, but let us be 
explicit: If $z_{j}$ is represented by 
$(v_{j},h_{j})$ for $j=1, 2$, then $z_{1}.z_{2}$ is represented by 
$(w,k)$, where $k$ is cartesian over $\pi (h_{1})$ and $w$ is 
vertical, and the square displayed commutes:
\begin{equation}\label{compx}\begin{diagram}\cdot &\rTo^{k}&\cdot & \rTo^{h_{2}}&\cdot\\
\dTo^{w}&&\dTo_{v_{2}}&&\\
\cdot&\rTo_{h_{1}}&\cdot &&\\
\dTo^{v_{1}}&&&&\\
\cdot
\end{diagram}\end{equation}
Such $k$ and $w$ exists (uniquely, up to unique vertical cartesian 
arrows): 
construct first $k$ as a cartesian lift of $\pi (h_{1})$, then use 
the universal property of cartesian arrows to construct $w$. (The 
square displayed will then actually be a pull-back diagram, thus the 
composition described will be the standard composition of spans.)

Composition of vh spans  does not give a definite vh span, 
but rather an 
equivalence class of vh spans. So referring to (\ref{compx}), the 
composite of  $\{(v_{1},h_{1})\}$ with the  
of $\{(v_{2},h_{2})\}$ is defined by
$$   \{(v_{1},h_{1})\}.\{(v_{2},h_{2})\}:= \{(w.v_{1},k.h_{2})\}.$$

There is a functor $\pi^{*}$ from $\X^{*}$ to $\B$; on objects, it agrees with $\pi 
: \X \to \B$; and   $\pi ^{*}(\{(v,h)\})= 
\pi (h)$. Note that if 
$v:X'\to X$ is vertical, the vh span $(v,1)$ represents a morphism 
$X \to X'$ in $\X^{*}$.  

Clearly, a vertical arrow in $\X^{*}$ has a unique representative 
span of the form $(v,1)$.
So  the fibres of $\pi^{*}:\X^{*} \to \B$ are 
canonically isomorphic
%\footnote{but not {\em equal} to; see the 
%Example at the end of the Section.} 
to the 
duals of the fibres of $\pi: \X \to \B$, i.e.\ $(\X ^{*})_{A} \cong (\X 
_{A})^{op}$; so 
 $\X^{*}$ is ``fibrewise dual'' to $\X$ (but is not 
in general dual to $\X$, since the functor $\pi^{*}:\X^{*}\to \B$ is 
still a covariant functor). 
The arrows in $\X^{*}$, we call {\em comorphisms}; it is ususally 
harmless to use the name ``comorphism'' also for a representing vh span $(v,h)$.

There are two special classes of comorphisms: the first class 
consists of those comorphisms that 
can be represented by a pair $(v,1)$ where $1$ is the relevant identity arrow. They 
are precisely the vertical arrows for $\X^{*}\to \B$.  -- The second class 
consists of those comorphisms that 
can be represented
  by a pair $(1,h)$ where $1$ is 
the relevant identity arrow. We shall see that these are precisely 
the cartesian morphisms in $\X^{*}$. '

We first note that if $(v,h)$ represents an arbitrary arrow in $\X^{*}$, 
then 
\begin{equation}\label{factx}(v,h) \in \{(v,1)\}.\{(1,h)\} ;
\end{equation}
this is witnessed by the diagram
$$\begin{diagram}\cdot& \rTo ^{1}&\cdot & \rTo^{h}&\cdot\\
\dTo^{1}&&\dTo_{1}&&\\
\cdot&\rTo_{1}&\cdot&&\\
\dTo^{v}&&&&\\
\cdot&&&&
\end{diagram}$$
since the upper left square is of the form considered in 
(\ref{compx}).

\begin{prop} An arrow $g$ is cartesian in $\X^{*}$ iff it admits a vh 
representative of the form $(1,h)$.
\end{prop}
{\bf Proof.} In one direction, let $(1,h)$ represent a comorphism 
$Y\to Z$ over $\beta \in \B$, and let 
$(v, k)$ represent  a comorphism $X\to Z$  over $\alpha. \beta$. We 
display these data as the full arrows in the following display (in 
$\X$ and $\B$):
$$ \begin{diagram}\cdot&&&&\\
\dTo^{v}&\rdDotsto_{k'} \rdTo(4,2)^{k}&&&\\
X&&Y & \rTo_{h}& Z\\
:&&: && :\\
\cdot&\rTo_{\alpha }&\cdot &\rTo _{\beta }&\cdot 
\end{diagram};$$
The dotted arrow  $k'$ comes about by using the universal property of the 
cartesian arrow $h$ in $\X$. Since $k$ and $h$ are Cartesian, then so 
is $k'$, by the Lemma \ref{lemmax}. So $(v,k')$ is a comorphism 
over $\alpha$, and $(v,k').(1.h) \equiv (v,k)$, and using the cancellation 
property of Cartesian arrows,  $(v,k')$ is easily 
seen to be the unique comorphism over $\alpha . \beta$ composing with 
$(1,h)$ to give  $(v,k)$.  

In the other direction, let $g$ be a cartesian arrow in $\X^{*}$. 
Let $(w,k)$ be an arbitrary representative of $g$. Then 
by (\ref{factx}), $g=\{(w,1)\}.\{(1,k)\}$. Since $g$ is assumed cartesian in 
$\X^{*}$, and $\{(1,k)\}$ is cartesian  by what is already proved, it 
follows from Lemma \ref{lemmax} that $\{(w,1)\}$ is 
cartesian. Since it is also vertical, it follows that it is an 
isomorphism in $\X^{*}$, hence $w$ is an isomorphism in $\X$. Since 
$k$ is cartesian in $\X$, $w^{-1}.k$ is cartesian as well, and
$$(w,k)\equiv (1,w^{-1}.k),$$
so $g$ has a representative of the claimed form.

\medskip

\begin{prop}The functor $\pi^{*}:\X^{*}\to \B$ is a fibration over $\B$
\end{prop}
{\bf Proof.} Let $\beta :A \to B$ be an arrow in $\B$, and let $Y$ be an 
object in over $B$. Since $\X \to \B$ is a fibration, there exists in 
$\X$ a cartesian arrow $h$ over $\beta$, and then the vh span $(1,h)$ 
represents, by the above, a cartesian arrow in $\X^{*}$ over $\beta$.

\medskip

Since $\X^{*}\to \B$ is a fibration, we may ask for its fibrewise 
dual $\X ^{**}$: 
\begin{prop}There is a canonical isomorphism over $\B$ between $\X$ 
and $\X^{**}$.
\end{prop}
{\bf Proof.} We describe an explicit functor $y: \X \to \X^{**}$.
Let us denote arrows in $\X^{*}$ by dotted arrows; they may be 
presented by vh spans $(v,h)$ in $\X$.
We first describe $y$ on vertical and cartesian arrows separately. For a 
vertical $v$ in $\X$, say $v:X'\to X$, we have the  
vh span $(v,1)$ in $\X$, which represents     a 
vertical 
arrow $\overline{v}: X'\dasharrow X$ in $\X^{*}$; thus we have a 
vh span $(\overline{v},1)$ in $\X^{*}$, which in turn represents
 a vertical arrow $X\to X'$ in $\X^{**}$. This arrow, we take as $y(v) \in \X^{**}$. 
Briefly, $y(v)= ((v,1),1)$. -- For a cartesian $h:X'\to Y$ (over 
$\beta$, say), we have a vh span $(1, h)$ in $\X$, which represents  a horizontal  arrow $\overline{h}: X'\dasharrow Y$ in $\X^{*}$ (cartesian 
over $\beta$); 
thus we have a vh span $(1,\overline{h})$ in $\X^{*}$, hence an 
arrow in $\X^{**}$, from $X'$ to $Y$ which we take as $y(h) \in \X^{**}$;
 briefly, $y(h)= (1,(1,v))$. 

Then, for a general $f:X\to Y$ in $\X$, we factor it $v.h$ with $v$ 
vertical and $h$ cartesian, and put $y(f):= y(v).y(h)$. We leave to 
the reader to verify that a different choice of $v$ and $h$ gives an equivalent 
vh span in $\X^{*}$, thus the same arrow in $\X^{**}$. 

Conversely, given an arrow $g: X \to Y$ in $\X^{**}$, represent it by 
a vh span  
in $\X^{*}$, $(\overline{v}, \overline{h})$,
$$\begin{diagram}X'& \rDotsto ^{\overline{h}}&Y\\
\dDotsto ^{\overline{v}}&&\\
X&&
\end{diagram}$$
Since $\overline{v}$ is vertical, we may pick a 
representative of $\overline{v}$ in the form $(v,1)$ with $v:X\to X'$,
 and since $\overline{h}$ is cartesian in $\X^{*}$, 
we may pick a representative of it if the form $(1,h)$, with $h:X'\to Y$ 
in $\X$. Then the composite $v.h : X\to Y$ makes sense in $\X$, and it goes by 
$y$ to the given $g$.

\medskip

\noindent{\bf Example.} Consider a group homomorphism $\pi : \X \to 
\B$.  It is a fibration iff $\pi$ is surjective. Assume 
this. Then the fibre (over the unique object $*$ of $\B$) is the 
kernel $\K$ of $\pi$. Every $h\in \X$ is Cartesian; the vertical 
arrows are those of $\K$. Then $\X^{*}$ is canonically isomorphic to 
$\X$. 
For, an element (arrow) $(v,h)$ of $\X^{*}$ may be presented by
either $(1, v^{-1}.h)$, so may be presented in 
the form $(1,x)$. 
The map $(v,h) \mapsto v^{-1}.h$ gives a canonical isomorphism  $J: \X^{*} \to 
\X$. This isomorphism preserves $\pi$; note that the $\pi$ for 
$\X^{*}$ takes $(v,h)$ to $\pi (h)$. Let us for clarity denote it 
$\pi'$, so $\pi' \{(v,h)\} = \pi (h)$. The kernel $\K'$ for $\pi'$ consists of elements which may be 
represented in the form $(v,1)$ with $v\in \K$, so $\K'$ may, as a 
set,  be 
identified with $\K$ by identifying $(v,1)\in \K' \subseteq \X^{*}$ with $v\in \K 
\subseteq \X$.
 But this identification is an anti-isomorphism, since 
$(v,1)$ by $J$ goes to $v^{-1}.1=v^{-1}$. So $\K '$ is identified as 
a group with $\K ^{op}$. Thus we have a diagram of 
group homomorphisms
$$\begin{diagram}\K^{op}&\rTo^{(-)^{-1}}_{\cong} &\K\\
\dTo^{i} &&\dTo_{\subseteq}\\
\X^{*}&\rTo^{J}_{\cong}&\X\\
\dTo^{\pi'}&&\dTo_{\pi}\\
\B & \rTo_{id}&\B
\end{diagram}$$ 
where $i(v)=\{(v,1)\}$. 
In case where $\B=1$, and $\X$ is the group $G$, the four maps of the 
top square are more explicitly`the four group isomorphisms
$$\begin{diagram}G^{op}&\rTo ^{(-)^{-1}_{i}}&G\\
\dTo^{v\mapsto {\{(v,1)\}}}&&\dTo_{=}\\
G^{*}&\rTo_{\{(v,h)\}\mapsto v^{-1}.h}^{J}&G
\end{diagram}$$
where the inverse of $J$ is given by $h\mapsto \{(1,h)\}$. If we 
denote the inverse of $J$ by $j$, we can write the information in 
this diagram more symmetrically:
$$\begin{diagram}G^{op}& \rTo^{i}&G^{*}&\lTo^{j}&G
\end{diagram}$$
with $i(v):= \{(v,1)\}$ and $j(h):= \{(1,h)\}$.
%Note that $G^{*}$ the here is the orbit of $G\times G$ under the 
%-- The reader should in particular 
%contemplate the special case where $\B$ is the one-element group.

\section{The case of a (pseudo-) functor $\B ^{op}\to Cat$}
 It is well known that a pseudofunctor $F: \B^{op} \to Cat$,
 gives rise to a fibration 
over $\B$. It is described in, say, \cite{Johnstone} B.1.3, or in \cite{Borceux}
 II.8.3. This fibration is known as the Grothendieck 
construction for $F$. We descibe it briefly in terms of the factorization 
system alluded to in Section 1.

Given a functor (or just a pseudo-functor) $F: \B ^{op} \to Cat$.
 Then we have a category $\X$ whose objects are pairs $(X,A)$ 
with $A$ an object of $\B$ and $X$ an object in $F(A)$. Arrows $(X,A) 
\to (Y,B)$ are pairs $(v, \alpha)$, where $\alpha :A \to B$ and $v: 
 X\to \alpha^{*}(Y)$ in $F(A)$ (and where $\alpha^{*}$ denotes the functor 
$F(\alpha): F(B) \to F(A)$). The functor $\pi : \X \to \B$ takes this 
arrow to $\alpha$.

Let us denote the arrow $(1_{\alpha^{*}(Y)}, \alpha )$ by $\alpha 
\dcoc Y$, thus
$$\begin{diagram}\alpha ^{*}(Y)&\rTo ^{\alpha \dcoc Y}&Y\end{diagram}$$
This is a Cartesian arrow over $\alpha$  in $\X$, and every Cartesian 
arrow is of this form modulo  unique vertical isomorphisms. There is 
then a canonical factorization of general arrows in $\X$, namely, the 
arrow given by a pair $(v, \alpha )$, as above, factors
as
$$\begin{diagram}(X,A)&&\\
\dTo^{(v, 1_{A})} &&\\ 
(\alpha^{*}(Y),A)&\rTo_{\alpha \dcoc Y}& (Y,B)
\end {diagram}.$$
Let $F'$ be $F$ followed by the dualization functor $Cat \to Cat$. 
Then a morphism over $\alpha$ in the fibration corresponding to $F'$, from $(X,A) $ to 
$(Y,B)$, is given similarly, but now with $v: \alpha^{*}(Y) \to X$, 
which in terms of the category $F(A)$ rather than $(F(A))^{op}$ may 
be displayed in terms of the vh span
$$\begin{diagram}(\alpha^{*}(Y),A)&\rTo^{\alpha \dcoc Y}& (Y,B)\\
\dTo^{(v,1_{A})}&&\\
(X,A)&&
\end{diagram},$$
and from this, it is clear that the fibration corresponding to $F'$ 
is isomorphic to $\X^{*}$ as described in the previous Sections.

\medskip

\small One motivation for the present note is to extract the pure 
category theory behind ``fibrewise 
contravariant functors" (like fibrewise duality for vector bundles), 
and ``star-bundle functors'', as in [Kolar et al, 1993] 41.2. This is 
still an ongoing project. 

\medskip

I cannot imagine that the constructions in the present note are not 
known, but I do not presently know of any available account.

\bigskip

\noindent
University of Aarhus, January 2015

\noindent mail: kock (at) math.au.dk

\end{document}